\title{Exact solution of nonlinear ordinary differential equations}
\author{
        Ming Tian Xu \\
                Department of Engineering Mechanics\\
        Shandong University\\
        Jinan City, 250014, \underline{P. R. China}
}
\date{\today}

\documentclass[12pt]{article}

\begin{document}
\maketitle

\begin{abstract}
At present, only some special differential equations have explicit analytical solutions. In general, no one thinks that it is possible to analytically find the exact solution of nonlinear equations. In this article based on the idea that the numerical scheme with zero truncation error can give rise to exact solution, a general formula for the exact solution of the initial value problem of nonlinear ordinary differential equations is obtained. Meanwhile, this formula enables us to construct a numerical scheme with zero truncation error for solving the nonlinear differential equations.
\end{abstract}
\section{Introduction}
At the very beginning, the research of differential equations was mainly devoted to develop various techniques for integrating particular types of differential equations and to express the solutions with elementary function or transcendental function. This upsurge was interrupted by Liouville's proof in 1841 that the solution of the special Riccati equation cannot be expressed via integration of elementary functions\cite{roge}. Since then the qualitative method, stability analysis method and dynamic system theory have become the mail tools to study the nonlinear differential equations\cite{roge,lyap,arno,hirs}. In addition, perturbation and singular perturbation method as well as transform group theory have been developed rapidly in the 20th century for investigating the nonlinear differential equations\cite{nayf,john,blum,ibra}. After the 1960s and 1970s, numerical methods for the solution of initial value problems of ordinary differential equations made enormous progress\cite{butc}. It is generally believed that numerical methods can only lead to numerical solutions with some errors. However, according to the definition of truncation error, when the truncation error of a numerical scheme vanishes, it should yield an exact solution. Based on this idea we attempt to obtain exact solutions of initial value problems of nonlinear differential equations.

\section{Exact solution of ordinary differential equations}
Consider the following system of ordinary differential equations:
\begin{eqnarray}
\frac{d{\mathbf{y}}}{dx}={\mathbf{f}}(x;{\mathbf{y}}) \label{e0}
\end{eqnarray}
with the initial condition:
\begin{eqnarray}
{\mathbf{y}}|_{x=x_0}={\mathbf{y}}_0\label{e1}
\end{eqnarray}
where ${\mathbf{y}}=(y_1,y_2,\cdots,y_n)^T$, $y_i(x) (i=1,2,\cdots,n)$ are unknown functions of ordinary differential equation system (\ref{e0}), $\mathbf{f}=(f_1,f_2,\cdots,f_n)^T$, $f_i(x;y_1,y_2,\cdots,y_n)$\\$(i=1,2,\cdots,n)$ are given smooth functions, ${\mathbf{y}}_0=(y_{10},y_{20},\cdots,y_{n0})^T$, $y_{i0}(i=1,2,\cdots,n)$ are known values. For the initial value problem described by Eqs.(\ref{e0}) and (\ref{e1}) we have the following Theorem.\\
{\textbf{Theorem 1}}:\\
If
\begin{eqnarray}
f_k(x;z_1,z_2,\cdots,z_n)-f_k(x_0;y_{10},y_{20},\cdots,y_{n0})\neq 0\label{ee0}
\end{eqnarray}
with
\begin{eqnarray}
z_k=y_{k0}+(x-x_0)f_k(x_0;{\mathbf{y}}_0) (k=1,2,\cdots,n)\label{e3}
\end{eqnarray}
The exact solution of the initial value problem (\ref{e0})-(\ref{e1}) is expressed as
\begin{eqnarray}
y_k(x)=y_{k0}+(x-x_0)[\alpha_kf_k(x_0;{\mathbf{y}}_0)+(1-\alpha_k)f_k(x;z_1,z_2,\cdots,z_n)]\label{ee1}
\end{eqnarray}
with\\
\begin{eqnarray}
\alpha_k=1-\frac{\sum_{m=2}^{\infty}\frac{(x-x_0)^{m-1}}{m!}D^{m-1}f_k(x;y_{10},y_{20},\cdots,y_{n0})}{f_k(x;z_1,z_2,\cdots,z_n)-f_k(x_0;y_{10},y_{20},\cdots,y_{n0})}\label{e2}\\
(k=1,2,\cdots,n)\nonumber
\end{eqnarray}
where
\begin{eqnarray}
D=\frac{\partial}{\partial x}+f_1\frac{\partial}{\partial y_1}+f_2\frac{\partial}{\partial y_2}+\cdots+f_n\frac{\partial}{\partial y_n}\label{e2a}
\end{eqnarray}
\emph{Proof}: Equivalently, Eq.(\ref{e0}) is rewritten into:
\begin{eqnarray}
\frac{dy_k}{dx}=f_k(x;y_1,y_2,\cdots,y_n)(k=1,2,\cdots,n)\label{e4}
\end{eqnarray}
Integrating both sides of the above equation yields
\begin{eqnarray}
y_k(x)=y_{k0}+\int_{x_0}^{x}f_k(x;y_1,y_2,\cdots,y_n)dx\label{e5}\\
(k=1,2,\cdots,n)\nonumber
\end{eqnarray}
The above integral equations can be discretized as
\begin{eqnarray}
\tilde{y}_k(x)&=&y_{k0}+(x-x_0)[\alpha_kf_k(x_0;{\mathbf{y}}_0)+(1-\alpha_k)f_k(x;z_1,z_2,\cdots,z_n)]\label{e6}\\
&&(k=1,2,\cdots,n)\nonumber
\end{eqnarray}
where $\alpha_k (k=1,2,\cdots,n)$ are undetermined parameters, $\tilde{y}_k(x)$ is the approximate function of $y_k(x)$. Note that when $\alpha_k(k=1,2,\cdots,n)$ take the value $0.5$, the algorithm (\ref{e6}) is nothing more than the second order Runge-Kutta method. Normally, even the values of the unknown functions at $x_0$ are exact, the values of $\tilde{y}_k(x) (k=1,2,\cdots,n)$ at the point $x$ have some errors in comparison with the exact solution $y_k(x)$, which is induced by the truncation error of the discretization scheme (\ref{e6}). Therefore, if one wants to get the exact solution from Eq.(\ref{e6}), then the truncation error of the algorithm (\ref{e6}) must vanish, that is,
\begin{eqnarray}
y_k(x)-y_{k0}-(x-x_0)[\alpha_kf_k(x_0;{\mathbf{y}}_{0})+(1-\alpha_k)f_k(x;z_1,z_2,\cdots,z_n)]=0\label{e7}
\end{eqnarray}
If $f_k(x;z_1,z_2,\cdots,z_n)-f_k(x_0;y_{10},y_{20},\cdots,y_{n0})\neq 0$, solving the above equation gives
\begin{eqnarray}
\alpha_k=\frac{-y_k(x)+y_{k0}+(x-x_0)f_k(x;z_1,z_2,\cdots,z_n)}{(x-x_0)[f_k(x;z_1,z_2,\cdots,z_n)-f_k(x_0;{\mathbf{y}}_0)]}\label{e8}
\end{eqnarray}
By expanding $y_k(x)$ to Taylor's series at the point $x_0$, we obtain
\begin{eqnarray}
y_k(x)&=&y_{k0}+(x-x_0)\frac{\partial y_k}{\partial x}|_{x=x_0}+\frac{(x-x_0)^2}{2!}\frac{\partial^2 y_k}{\partial x^2}|_{x=x_0}\nonumber\\
&&+\cdots+\frac{(x-x_0)^m}{m!}\frac{\partial^m y_k}{\partial x^m}|_{x=x_0}+\cdots\label{e9}
\end{eqnarray}
By using Eq.(\ref{e4}), we have
\begin{eqnarray}
\frac{\partial^m y_k}{\partial x^m}=D^{m-1}f_k (m=1,2,\cdots\cdots)\label{e10}
\end{eqnarray}
where
\begin{eqnarray}
D=\frac{\partial}{\partial x}+f_1\frac{\partial}{\partial y_1}+f_2\frac{\partial}{\partial y_2}+\cdots+f_n\frac{\partial}{\partial y_n}\label{e10a}
\end{eqnarray}
Substituting Eq.(\ref{e10}) into (\ref{e9}) yields
\begin{eqnarray}
y_k(x)&=&y_{k0}+(x-x_0)f_k|_{x=x_0}+\frac{(x-x_0)^2}{2!}D f_k|_{x=x_0}\nonumber\\
&&+\cdots+\frac{(x-x_0)^m}{m!}D^{m-1} f_k|_{x=x_0}+\cdots\label{e11}
\end{eqnarray}
Inserting the above equation into Eq.(\ref{e8}) gives
\begin{eqnarray}
\alpha_k=1-\frac{\sum_{m=2}^{\infty}\frac{(x-x_0)^{m-1}}{m!}D^{m-1}f_k|_{x=x_0}}{f_k(x;z_1,z_2,\cdots,z_n)-f_k(x_0;{\mathbf{y}}_0)}\label{e12}
\end{eqnarray}
Therefore, when the parameters $\alpha_k (k=1,2,\cdots,n)$ are calculated by Eq.(\ref{e12}), expression (\ref{e6}) provides us the exact solution of the initial value problem (\ref{e0})-(\ref{e1}).\\
{\textbf{Remark 1}}. Note that the exact solution expressed by Eqs.(\ref{e6}) and (\ref{e12}) is applicable to the system of linear ordinary differential equations.\\
{\textbf{Remark 2}}. Since the following system of the $m$-th order differential equations 
\begin{eqnarray}
\frac{d^m{\mathbf{y}}}{dx^m}={\mathbf{g}}(x;{\mathbf{y}},\frac{d{\mathbf{y}}}{dx},\cdots,\frac{d^{m-1}{\mathbf{y}}}{dx^{m-1}})\label{e13}
\end{eqnarray} 
where ${\mathbf{g}}=(g_1,g_2,\cdots,g_n)^T$ is a smooth function, can be written into the system (\ref{e0}),
the exact solution expressed by Eqs.(\ref{e6}) and (\ref{e12}) can also be applied to solve the system (\ref{e13}).\\
\mbox{ }\mbox{ }Theorem 1 can allow us to construct a numerical scheme which yields the solution with the exact solution's accuracy.\\ 
{\textbf{Corollary 1}}. One can use the following algorithm to solve the initial value problem (\ref{e0}-\ref{e1})
\begin{eqnarray}
y_{k,i+1}=y_{k,i}+h[\alpha_k f_{k,i}+(1-\alpha_k)f_k(x_{i+1};z_{1,i},z_{2,i},\cdots,z_{n,i})]\label{e14}\\
 (k=1,2,\cdots,n)\nonumber
\end{eqnarray}
where $x_{i+1}=x_i+h$, $h$ is the time step, $y_{k,i}=y_k(x_i)$, $f_{k,i}=f_k(x_i;y_1(x_i),y_2(x_i),$
$\cdots,y_n(x_i))$, $z_{j,i}=z_j(x_i)$. From Theorem 1 one can see that the truncation error of scheme (\ref{e14}) vanishes, thus it can achieve the exact solution's accuracy.
\section{Examples}
{\textbf{Example 1.}} Consider the following initial value problem of linear differential equation:
\begin{eqnarray}
\biggl\{\begin{array}{c}
\frac{dy}{dx}=e^x\label{e15a}\\
y(0)=1\label{e15b}\\
\end{array}
\end{eqnarray}
Obviously, the exact solution of this problem is $y=e^x$. Next we shall examine the solution expressed by Eqs.(\ref{e6}) and (\ref{e12}). Firstly, 
\begin{eqnarray}
Df=(\frac{\partial}{\partial x}+e^x\frac{\partial}{\partial y})e^x=e^x\label{e16}
\end{eqnarray}
subsequently, $D^kf=e^x (k=1,2,\cdots\cdots)$. Substituting these results into Eq.(\ref{e12}) yields
\begin{eqnarray}
\alpha&=&1-\frac{\frac{x}{2!}+\frac{x^2}{3!}+\cdots+\frac{x^{k-1}}{k!}+\cdots}{e^x-1}\nonumber\\
&=&1-\frac{1}{x}+\frac{1}{e^{x}-1}\label{e17b}
\end{eqnarray}
Inserting the above into Eq.(\ref{e6}) gives
\begin{eqnarray}
y(x)&=&1+x[1-\frac{1}{x}+\frac{1}{e^{x}-1}+(\frac{1}{x}-\frac{1}{e^{x}-1})e^{x}]\nonumber\\
&=&e^x\label{e18}
\end{eqnarray}
where the initial condition in (\ref{e15b}) has been used. This shows that the solution expressed by Eqs.(\ref{e6}) and (\ref{e12}) is realy the exact solution of the initial value problem (\ref{e15a}).\\
{\textbf{Example 2.}} Consider the following initial value problem of nonlinear differential equation
\begin{eqnarray}
\biggl\{\begin{array}{c}
\frac{dy}{dx}=y^2\label{e19a}\\
y(0)=1\label{e19b}
\end{array}
\end{eqnarray}
The exact solution of this problem is $y=\frac{1}{1-x}$. Denote $f(x,y)=y^2$. Note that 
\begin{eqnarray}
Df=(\frac{\partial}{\partial x}+y^2\frac{\partial}{\partial y})y^2=2y^3\label{e20}
\end{eqnarray}
Subsequently, we have
\begin{eqnarray}
D^kf=(k+1)!y^{k+2}(k=1,2,\cdots\cdots)\label{e21}
\end{eqnarray}
Substituting Eq.(\ref{e21}) into (\ref{e12}) gives
\begin{eqnarray}
\alpha=1-\frac{x+x^2+\cdots+x^{k-1}+\cdots}{(1+x)^2-1}\label{e22}
\end{eqnarray}
Inserting Eq.(\ref{e22}) into (\ref{e6}) yields
\begin{eqnarray}
y(x)&=&1+x+x^2+x^3+\cdots+x^k+\cdots\nonumber\\
&=&\frac{1}{1-x}\label{e23}
\end{eqnarray}
Therefore, Eqs.(\ref{e6}) and (\ref{e12}) realy give rise to the exact solution of the initial value problem (\ref{e19a}).\\
{\textbf{Example 3.}} Consider the following Riccati equation:
\begin{eqnarray}
\frac{dy}{dx}=P(x)+Q(x)y+R(x)y^2\label{e24}
\end{eqnarray}
Firstly, by applying the transform\cite{ibra1}
\begin{eqnarray}
z=-Ry-0.5(\frac{R'}{R}+Q)\nonumber
\end{eqnarray}
where $R'$ is the first order derivative of $R$. Eq.(\ref{e24}) becomes
\begin{eqnarray}
\frac{dz}{dx}+z^2=0.25(\frac{R'}{R}+Q)^2-0.5(\frac{R''}{R}-\frac{{R'}^2}{R^2}+Q')-PR\label{e25}
\end{eqnarray}
Set $P=e^x-e^{3x}$,$Q(x)=2e^{2x}$,$R=-e^x$, Eq.(\ref{e25}) becomes
\begin{eqnarray}
z'+z^2=0.25\label{e26}
\end{eqnarray}
and $z$ satisfies the initial condition $z(0)=0$, that is, $y(0)=1.5$. The exact solution of this problem is that $y=e^x+e^{-x}$. We use the following finite sum 
\begin{eqnarray}
\sum_{m=1}^{8}\frac{x^m}{(m+1)!}D^mf|_{x=0}\label{e26a}
\end{eqnarray}
with
\begin{eqnarray}
Df&=&-0.5z+2z^3\nonumber\\
D^2f&=&-\frac{1}{8}+2z^2-3!z^4\nonumber\\
D^3f&=&z-10z^3+4!z^5\nonumber\\
D^4f&=&\frac{1}{4}-\frac{17}{2}z^2+60z^4-5!z^6\nonumber\\
D^5f&=&-\frac{17}{4}z+77z^3-420z^5+6!z^7\nonumber\\
D^6f&=&-\frac{17}{16}+62z^2-756z^4+3360z^6-7!z^8\nonumber\\
D^7f&=&31z-880z^3+8064z^5-30240z^7+8!z^9\nonumber\\
D^8f&=&\frac{31}{4}-691z^2+12720z^4-93240z^6+302400z^8-9!z^{10}\nonumber
\end{eqnarray}
to approximate the series in Eq. (\ref{e12}) for calculating the weighting parameter $\alpha$. The scheme (\ref{e14}) is employed to calculate the solution with the step $h=0.1$. The obtained solution in this way and the exact solution is listed in Table 1. One can see that although the finite sum is used to approximate the series in Eq.(\ref{e12}) the obtained solution still achieves the exact solution's accuracy. If from the initial condition $y(0)=1.5$ we directly compute $y(1)$ by Eq.(\ref{e6}), the numerical result is $y(1)=2.98722397263779$, the exact solution is $y(1)=2.98722324982904$, the relative error is about $2.95E-7$, therefore, a quite high accuracy is also achieved. This shows that the power series in Eq.(\ref{e12}) has a fast convergence rate.\\
\begin{table}
\caption{Comparison of the solution obtained by the present method and exact solution for Example 3}
\label{tab2}       
\begin{center}
\begin{tabular}{llll}
\hline\noalign{\smallskip}
$x$ & Present method & Exact solution \\
\noalign{\smallskip}\hline\noalign{\smallskip}
$0.1$ & $1.58019173059671$ & $1.58019173059671$\\
$0.2$ & $1.67156876084769$ & $1.67156876084769$\\
$0.3$ & $1.77541629076434$ & $1.77541629076434$\\
$0.4$ & $1.89313703752882$ & $1.89313703752882$\\
$0.5$ & $2.02626193949827$ & $2.02626193949827$ \\
$0.6$ & $2.17646249416471$ & $2.17646249416471$ \\
$0.7$ & $2.34556493530231$ & $2.34556493530231$ \\
$0.8$ & $2.53556644736486$ & $2.53556644736486$ \\
$0.9$ & $2.74865360853195$ & $2.74865360853195$ \\
$1.0$ & $2.98722324982904$ & $2.98722324982904$ \\
\noalign{\smallskip}\hline
\end{tabular}
\end{center}
\end{table}
{\textbf{Example 4}} Consider the following inital value problem of the second order linear differential equation:
\begin{equation}{\label{e27}}
\biggl\{\begin{array}{c}
\frac{d^2\phi}{dx^2}+\omega^2\phi=0\\
x=0;\phi=\phi_0,\frac{d\phi}{dx}=0\\
\end{array}
\end{equation}
where $\omega$ is a system parameter.
Set $y_1=\phi$, $y_2=d\phi/dt$, Eq.(\ref{e27}) can be rewritten into:
\begin{equation}{\label{e28}}
\biggl\{\begin{array}{c}
\frac{dy_1}{dx}=f_1(x;y_1,y_2)\\
\frac{dy_2}{dx}=f_2(x;y_1,y_2)\\
\end{array}
\end{equation} 
with
\begin{eqnarray}
f_1(x;y_1,y_2)&=&y_2\label{e29a}\\
f_2(x;y_1,y_2)&=&-\omega^2y_1\label{e29b}
\end{eqnarray}
From Eqs.(\ref{e29a}) and (\ref{e29b}), we have
\begin{eqnarray}
Df_1=(\frac{\partial}{\partial x}+f_1\frac{\partial}{\partial y_1}+f_2\frac{\partial}{\partial y_2})f_1=-\omega^2y_1\label{e30a}\\
D^2f_1=(\frac{\partial}{\partial x}+f_1\frac{\partial}{\partial y_1}+f_2\frac{\partial}{\partial y_2})(-\omega^2 y_1)=-\omega^2y_2\label{e30b}
\end{eqnarray}
Subsequently, we obtain
\begin{eqnarray}
D^{2k-1}f_1=(-1)^k\omega^{2k}y_1\label{e31a}\\
D^{2k}f_1=(-1)^k\omega^{2k}y_2\label{e31b}\\
(k=1,2,\cdots)\nonumber
\end{eqnarray}
Similarly we get
\begin{eqnarray}
D^{2k-1}f_2=(-1)^k\omega^{2k}y_2\label{e32a}\\
D^{2k}f_2=(-1)^{k+1}\omega^{2k+2}y_1\label{e32b}\\
(k=1,2,\cdots)\nonumber
\end{eqnarray}
Substituting Eqs.(\ref{e31a})-(\ref{e32b}) into (\ref{e12}) yields
\begin{eqnarray}
\alpha_1=1+\frac{\sum_{k=1}^{\infty}[(-1)^k\frac{x^{2k-1}\omega^{2k}}{(2k)!}y_1(0)+(-1)^k\frac{(x\omega)^{2k}}{(2k+1)!}y_2(0)]}{x\omega^2y_1(0)}\nonumber\\
=1+\frac{\omega(cos(x\omega)-1)y_1(0)+(sin(x\omega)-x\omega)y_2(0)}{x^2\omega^3y_1(0)}\label{e33}
\end{eqnarray}
\begin{eqnarray}
\alpha_2=1+\frac{\sum_{k=1}^{\infty}[(-1)^k\frac{x^{2k-1}\omega^{2k}}{(2k)!}y_2(0)+(-1)^{k+1}\frac{x^{2k}\omega^{2k+2}}{(2k+1)!}y_1(0)]}{x\omega^2y_2(0)}\nonumber\\
=1+\frac{(cos(x\omega)-1)y_2(0)-\omega(sin(x\omega)-x\omega)y_1(0)}{x^2\omega^2y_2(0)}\label{e34}
\end{eqnarray}
Substituting Eqs.(\ref{e33}) and (\ref{e34}) into Eq.(\ref{e6}) gives
\begin{eqnarray}
y_1=\phi(x)=\phi_0cos(\omega x)\label{e35}
\end{eqnarray}
This is nothing more than the exact solution of the differential equation problem (\ref{e27}).\\
{\textbf{Example 5.}} Consider the following van der Pol equation
\begin{eqnarray}{\label{e36}}
\biggl\{\begin{array}{c}
\frac{d^2x}{dt^2}+2\mu (x^2-1)\frac{dx}{dt}+x=0\\
x(0)=1.8,\frac{dx}{dt}|_{t=0}=1.8\\
\end{array}
\end{eqnarray}
Set $y_1=x$, $y_2=dx/dt$. The above problem can be rewritten as follows
\begin{eqnarray}{\label{e37}}
\biggl\{\begin{array}{ccc}
\frac{dy_1}{dt}&=&y_2\\
\frac{dy_2}{dt}&=&2\mu(1-y_1^2)y_2-y_1\\
\end{array}
\end{eqnarray}
Therefore, $f_1=y_2$, $f_2=2\mu(1-y_1^2)y_2-y_1$. In accordance with Eq.(\ref{e2a}), we obtain
\begin{eqnarray}{\label{e38}}
Df_1&=&py_1-y_1\nonumber\\
D^2f_1&=&p^{(1)}y_2^2-y_2+p^2y_2-py_1\nonumber\\
D^3f_1&=&p^{(2)}y_2^3+4pp^{(1)}y_2^2+(-2p+p^3-3p^{(1)}y_1)y_2+(1-p^2)y_1\nonumber\\
D^4f_1&=&(4p^{(1)^2}+7pp^{(2)})y_2^3+(-5p^{(1)}+11p^2p^{(1)}-6p^{(2)}y_1)y_2^2\nonumber\\
&&+(-13pp^{(1)}y_1+1-3p^2+p^4)y_2+2py_1-p^3y_1+3p^{(1)}y_1^2\nonumber\\
D^5f_1&=&15p^{(1)}p^{(2)}y_2^4+(-11p^{(2)}+34pp^{(1)^2}+32p^2p^{(2)})y_2^3\nonumber\\
&&+(-25p^{(1)^2}y_1-46pp^{(2)}y_1-29pp^{(1)}+26p^3p^{(1)})y_2^2\nonumber\\
&&+(18p^{(1)}y_1+3p-38p^2p^{(1)}y_1-4p^3+15p^{(2)}y_1^2+p^5)y_2\nonumber\\
&&+13pp^{(1)}y_1^2-y_1+3p^2y_1-p^4y_1\nonumber\\
D^6f_1&=&15p^{(2)^2}y_2^5+(34p^{(1)^3}+192pp^{(1)}p^{(2)})y_2^4+(-156p^{(1)}p^{(2)}y_1\nonumber\\
&&-54p^{(1)^2}-108pp^{(2)}+180p^2p^{(1)^2}+122p^3p^{(2)})y_2^3+(81p^{(2)}y_1+21p^{(1)}\nonumber\\
&&-228pp^{(1)^2}y_1-226p^2p^{(2)}y_1-108p^2p^{(1)}+57p^4p^{(1)})y_2^2+(63p^{(1)^2}y_1^2\nonumber\\
&&+120pp^{(2)}y_1^2+108pp^{(1)}y_1+6p^2-94p^3p^{(1)}y_1-5p^4+p^6)y_2-18p^{(1)}y_1^2\nonumber\\
&&-3py_1+38p^2p^{(1)}y_1^2+4p^3y_1-15p^{(2)}y_1^3-p^5y_1\nonumber\\
D^7f_1&=&(294p^{(1)^2}p^{(2)}+267pp^{(2)^2})y_2^5+(-231p^{(2)^2}y_1-372p^{(1)}p^{(2)}\nonumber\\
&&+496pp^{(1)^3}+1494p^2p^{(1)}p^{(2)})y_2^4+(102p^{(2)}-364p^{(1)^3}y_1-2144pp^{(1)}p^{(2)}y_1\nonumber\\
&&-658p^2p^{(2)}+768p^3p^{(1)^2}-606pp^{(1)^2}+423p^4p^{(2)})y_2^3+(714p^{(1)}p^{(2)}y_1^2\nonumber\\
&&+396p^{(1)^2}y_1+834pp^{(2)}y_1+162pp^{(1)}-1398p^2p^{(1)^2}y_1-912p^3p^{(2)}y_1\nonumber\\
&&-330p^3p^{(1)}+120p^5p^{(1)})y_2^2+(-225p^{(2)}y_1^2-81p^{(1)}y_1-3p+595pp^{(1)^2}y_1^2\nonumber\\
&&+610p^2p^{(2)}y_1^2+412p^2p^{(1)}y_1+10p^3-213p^4p^{(1)}y_1-6p^5+p^7)y_2\nonumber\\
&&-63p^{(1)^2}y_1^3-120pp^{(2)}y_1^3-108pp^{(1)}y_1^2-6p^2y_1+94p^3p^{(1)}y_1^2\nonumber\\
&&+5p^4y_1-p^6y_1\nonumber
\end{eqnarray}
where
\begin{eqnarray}
p&=&2\mu(1-y_1^2)\nonumber\\
p^{(1)}&=&-4\mu y_1\nonumber\\
p^{(2)}&=&-4\mu\nonumber
\end{eqnarray}
And it is obvious that
\begin{eqnarray}
D^kf_2=D^{k+1}f_1 (k=1,2,\cdots)\label{e39}
\end{eqnarray}

\begin{table}
\caption{Comparison of the solution obtained by the present method and the fourth order Runge-Kutta mthod for Example 5}
\label{tab3}       
\begin{center}
\begin{tabular}{llll}
\hline\noalign{\smallskip}
$x$ & Present method & Runge-Kutta & Runge-Kutta \\
 &$h=0.1$&$h=1.0E-5$&$h=1.0E-7$\\
\noalign{\smallskip}\hline\noalign{\smallskip}
$0.1$ & $1.96652400267220$ & $1.96652401307289$ & $1.96652400271259$\\
$0.2$ & $2.10473374763406$ & $2.10473378177567$ & $2.10473374807376$\\
$0.3$ & $2.21351635642977$ & $2.21351641534912$ & $2.21351635650173$\\
$0.4$ & $2.29291198605060$ & $2.29291206256096$ & $2.29291198590522$\\
$0.5$ & $2.34390221872474$ & $2.34390229889988$ & $2.34390221783419$\\
$0.6$ & $2.36813285441244$ & $2.36813292036275$ & $2.36813285115724$\\
$0.7$ & $2.36763255340805$ & $2.36763258982445$ & $2.36763254899136$\\
$0.8$ & $2.34457029918107$ & $2.34457029110568$ & $2.34457029362396$\\
$0.9$ & $2.30107218816106$ & $2.30107212257853$ & $2.30107218096084$\\
$1.0$ & $2.23909995816732$ & $2.23909982596943$ & $2.23909994999267$\\
\noalign{\smallskip}\hline
\end{tabular}
\end{center}
\end{table}
Then we use the finite sums 
\begin{eqnarray}
\sum_{m=2}^8\frac{(x-1.8)^{m-1}}{m!}D^{m-1}f_1|_{x=1.8}\nonumber
\end{eqnarray} 
and 
\begin{eqnarray}
\sum_{m=2}^7\frac{(x-1.8)^{m-1}}{m!}D^{m-1}f_2|_{x=1.8}\nonumber
\end{eqnarray}
instead of the power series in Eq.(\ref{e12}) for calculating the weighting parameters $\alpha_1$ and $\alpha_2$, respectively. And scheme (\ref{e14}) with the step $h=0.1$ is employed to solve the initial value problem (\ref{e36}). The obtained results and the numerical solutions given by the fourth order Range-Kutta method are listed in Table {\ref{tab3}. One can see that scheme (\ref{e14}) has achieved a quite high accuracy. This shows that the power series in Eq.(\ref{e12}) has a very fast convergence rate. 
\section{Concluding remarks}
Based on the idea that a numerical scheme with zero truncation error can lead to exact solution of differential equations, an analytical expression of the exact solution of the initial value problem of nonlinear ordinary differential equations in the form of power series is derived by modifying the second order Runge Kutta scheme. The results show that the power series in the expression of the exact solution has a fast convergence rate.



\end{document}